\documentclass[12pt]{article}
\usepackage[final]{epsfig}	
\usepackage{amsfonts,amssymb,latexsym,amscd, theorem,epsfig,graphics}

\newtheorem{theorem}{Theorem}
\newtheorem{lemma}{Lemma}[section]
\newtheorem{remark}{Remark}[section]

\newtheorem{corollary}{Corollary}[section]
\newtheorem{conjecture}{Conjecture}[section]

\begin{document}
\newcommand{\eps}{{\varepsilon}}
\newcommand{\proofend}{$\Box$\bigskip}
\newcommand{\C}{{\mathbf C}}
\newcommand{\Q}{{\mathbf Q}}
\newcommand{\R}{{\mathbf R}}
\newcommand{\Z}{{\mathbf Z}}
\newcommand{\RP}{{\mathbf {RP}}}
\newcommand{\CP}{{\mathbf {CP}}}

\title {Existence and non-existence of skew branes}
\author{Serge Tabachnikov\thanks{Partially supported by NSF}\ \ and Yulia Tyurina\\
{\it Department of Mathematics, Penn State}\\
{\it University Park, PA 16802, USA}
}
\date{}
\maketitle

\section{Introduction}

In mid-1960s, H. Steinhaus conjectured that every closed smooth curve in 3-dimensional space has parallel tangent lines. Shortly after that, B. Segre \cite{Se 1, Se 2} constructed examples of  curves without parallel tangent lines but showed that no such curve can lie on  the unit sphere (and therefore, an ellipsoid). Call a curve without parallel tangents a {\em skew loop}. Geometrical and topological study of skew loops and their multi-dimensional analogs has become an active research subject.

We briefly survey the available results. The supply of skew loops is abundant: for example, every knot type can be realized by one \cite{Wu}. The aversion of skew loops to ellipsoids was extended to convex quadric surfaces in \cite{G-S} and to non-convex ones in \cite{Ta}, see also \cite{Gh2,S-S}. Ghomi and Solomon \cite{G-S}  prove a  converse statement: if a convex surface is not quadratic then it carries a skew loop. Another non-existence result \cite{Ta}:  no skew loop lies on a ruled developable disc.

A multi-dimensional version of a skew loop is called a {\em skew brane}.\footnote{The terms ``skew loop" and ``skew brane" were coined by M. Ghomi and B. Solomon; see  \cite{Gh1, G-T, S-T} for the study of other classes of non-degenerate embeddings of manifolds into affine and projective spaces}
 A skew brane $f: M^n \to \R^{n+2}$ is an immersion  such that the tangent spaces $df (T_{x} M)$ and $df (T_{y} M)$ are not parallel  for all $x \neq y$. 
Pairs of parallel tangent spaces correspond to self-intersections of the image of the tangent Gauss map $M \to G_{n} (n+2)$ in the Grassman manifold of $n$-dimensional subspaces in $\R^{n+2}$. Since dim  $G_{n} (n+2) = 2n$, these self-intersections generically occur in isolated points and cannot be destroyed by small perturbations of $f$. If $M$ is oriented, one distinguishes the cases when parallel tangent spaces have the same or the opposite orientations; we refer to the former as positively and the latter as negatively parallel tangent spaces.

The aversion of skew loops to quadric surfaces extends to skew branes: no skew brane can lie on a quadratic hypersurface of any signature \cite{Ta, S-S}; for spheres, this was proved in \cite{Wh}. The paper by Lai \cite{Lai2} has been overlooked in the more recent literature on the subject.  Unaware of the work by Segre, Lai provides an example of a skew loop in $\R^3$ and of a torus $T^2 \subset \R^4$, free from pairs of negatively parallel tangent spaces. The main result of \cite{Lai2} is the following  theorem on non-existence of skew branes.

\begin{theorem} \label{nonexist}
Let $M^{2n} \subset \R^{2n+2}$ be a closed oriented embedded submanifold  with non-zero Euler characteristic $\chi$. Then there exists a pair of distinct points $x,y \in M$ such that the tangent spaces $T_{x} M$ and $T_{y} M$ are negatively parallel. For a generic submanifold $M$, the number of such (unordered) pairs $(x,y)$ is not less than $\chi^2/4$.
\end{theorem}

We give a streamlined proof in Section \ref{pf1}. By a generic submanifold we mean the one whose Gauss map  is an immersion with transversal self-intersections; see Section \ref{pf1} for a justification. 

The papers by Lai \cite{Lai1,Lai2}  continue  the work by Blaschke \cite{Bl} and Chern and Spanier \cite{C-S}  concerning  embedded  surfaces in $\R^4$. In this case, we add to Theorem \ref{nonexist} the following result.  Let $M\subset \R^4$ be a closed oriented immersed surface of genus $g$. Assume that $M$ is generic in the sense that  self-intersection of $M$ are transversal. Then each double point is assigned a sign; denote by $d$ the algebraic number of double points.  

\begin{theorem} \label{nonexist2}
$M$ has at least $|d^2-(1-g)^2|$  pairs of negatively parallel tangent planes.
\end{theorem}

In Section \ref{pf2} we construct examples of skew branes. Let $\R^{2n}$ be a subspace in $\R^{2n+1}$ and  $S^{2n-1} \subset  \R^{2n}$ the unit sphere. We consider this sphere as a codimension 2 submanifold in $2n+1$-dimensional space.

\begin{theorem} \label{exist}
There exists a small perturbation of $S^{2n-1}$ in $\R^{2n+1}$, free from pairs of parallel tangent spaces. 
\end{theorem}

For $n=1$, a much stronger result is proved in \cite{G-S}: given a non-centrally symmetric smooth closed curve $\gamma$ in the horizontal plane, there exists a skew loop on the cylinder over $\gamma$ that projects diffeomorphically on $\gamma$ (``cylinder lemma"). We believe that an analog of the cylinder lemma holds for odd-dimensional spheres; see Conjecture \ref{cylinder} and a much more general Conjecture \ref{oneform}.

We also construct a skew torus in 4-dimensional space (our construction is different from that in \cite{Lai2}). Let $T_0^2\subset \R^4$ be the standard torus which is the product of two unit circles in the plane. 

\begin{theorem} \label{torus}
There exists a small perturbation of $T_0^2$ in $\R^4$, free from pairs of parallel tangent planes. 
\end{theorem}

Finally, we construct  an immersed 2-dimensional sphere in $\R^4$ with one double point and no pairs of negatively parallel tangent planes; this shows that Theorem \ref{nonexist2} is sharp, at least in the spherical case. Let $M_0$ be an immersed sphere in $\R^4$ given by the equation 
$$
(\alpha,h) \mapsto (1-h^2) (\cos \alpha, \sin \alpha, h \cos \alpha, h \sin \alpha)
$$
where $(\alpha,h)$ are the cylindrical coordinates on the unit sphere; the origin is the double point, the image of both poles given by $h=\pm1$. 

\begin{theorem} \label{2sph}
There exists a small perturbation of $M_0$, free from pairs of negatively parallel tangent planes. 
\end{theorem}

\smallskip

{\bf Acknowledgments}. We are grateful to D. Fuchs, M. Ghomi and B. Solomon for their interest and help.

\section{Topological obstructions to the existence of skew branes} \label{pf1}

In this section we prove Theorems \ref{nonexist} and \ref{nonexist2}. The proof of the former in \cite{Lai2} is based on a detailed analysis of the Schubert cell decomposition of  Grassmanian manifolds made in \cite{Lai1}; our proofs  use characteristic classes of vector bundles and are more straightforward. 

Consider the tangent Gauss map $F: M^{2n} \to G_{2n}^+ (2n+2)$ to the Grassman manifold of oriented $2n$-dimensional subspaces in $\R^{2n+2}$. Let $\sigma$ be the involution of the Grassmanian that inverts the orientation of every $2n$-dimensional subspace. The proofs consist of computing the homology class $F_*[M] \in H_{2n} (G_{2n}^+ (2n+2),\Z)$ where $[M]$ is the fundamental class of $M$ and of the homology intersection number $F_*[M]  \cap \sigma_* F_*[M]$. We suppress the coefficients from the notation of homology and cohomology.

Let $p$ and $q$ be even, and consider the Grassmanian $G_p^+(p+q)$. Denote by 
$\xi^p$ and $\nu^q$ the tautological vector bundles: the former has the oriented $p$-plane $E$, and the latter its orthogonal complement, as the fiber over $E$, considered as a point in $G_p^+(p+q)$. The bundles are oriented and so is the Grassmanian. The  tangent bundle is expressed as follows:  $T  G_p^+(p+q) = {\rm Hom}(\xi,\nu)$. Denote the fundamental class  by $c \in H_{pq} (G_p^+(p+q))$ and let $x=e(\xi)\in H_{p} (G_p^+(p+q))$ and $y=e(\nu)\in H_{q} (G_p^+(p+q))$ be the Euler classes. Denote by dot the pairing between homology and cohomology classes.

\begin{lemma} \label{multable}
One has: $c \cdot x^q=2,\ c \cdot y^p = 2 (-1)^{pq/4},\ xy=0$.
\end{lemma}

\noindent  {\bf Proof}. The homological Euler class of an oriented vector bundle $\eta$ over an oriented manifold $M$ can be found as follows: choose a generic section of $\eta$ and let $N\subset M$ be the set of  zeros of this section. Then $N$ is a submanifold,  oriented as the intersection of this section and the zero section.
The Euler class $e(\eta)$ is dual to the homology class $[N]$.

Fix a vector $e\in R^{p+q}$ and project it to each $p$-plane. This gives a section of $\xi$ whose zero manifold $N_e$ consists of $p$-planes orthogonal to $e$; one has: $N=G_p^+(p+q-1)$. The class $[N] \cap \dots \cap [N]$ ($q-1$ times) is dual to $x^{q-1}$; it is 
represented by the set $Q$ consisting of $p$-planes that are orthogonal to fixed $q-1$ vectors. Hence $Q= G_p^+ (p+1) = S^p$. One has: $c \cdot x^q = [Q] \cdot x$, 
which is the Euler number of the restriction of $\xi$ on $Q$. The latter is the tangent bundle of the sphere, and its Euler number is 2. 

Consider the map $d: G_p^+(p+q) \to G_q^+(p+q)$ that takes an oriented  $p$-plane to its orthogonal complement. Denote by $\xi'$ and $\nu'$ the tautological bundles over $G_q^+(p+q)$ and $x',y'$ their Euler classes. Then $d^*(\xi')=\nu, d^*(\nu')=\xi$ and $d^*(x')=y, d^*(y')=x$. Let $c'$ be the fundamental class of $G_q^+(p+q)$.
\smallskip

\noindent {\bf Claim}:  $d_*(c)=(-1)^{pq/4} c'$.
\smallskip

 Assuming this claim and interchanging $p$ and $q$, one has:
$$
2=c' \cdot (x')^p = (-1)^{pq/4} d_*(c)\cdot (x')^p = (-1)^{pq/4} c\cdot d^*(x')^p = (-1)^{pq/4} c \cdot y^p,
$$
as stated.

To prove that $xy=0$, fix a vector $f \in R^{p+q}$ and, for every $p$-plane $E$, project $f$ to $E^{\perp}$. This gives a section of $\nu$ whose zero manifold $P_f$ consists  of $p$-planes that contain $f$. The homology class, dual to $xy$, is  $[N_e \cap P_f]$. If $f$ is not orthogonal to $e$ then this intersection is empty, and hence $xy=0$.

It remains to prove the above claim. One has: 
$$
d^*(T G_q^+(p+q)) = d^*({\rm Hom}(\xi',\nu')) = {\rm Hom} (\nu,\xi).
$$
The bundles ${\rm Hom} (\xi,\nu)$ and ${\rm Hom} (\nu,\xi)$ are naturally isomorphic, and we need to find whether this isomorphism preserves orientations. 

The orientation of the former bundle is determined by the identification of a $p\times q$ matrix with the $pq$-vector formed by its first, second, etc., rows. Likewise, for the latter bundle, one considers the $pq$-vector formed by its first, second, etc., columns. A matrix element $(i,j)$ has positions $(i-1)p+j$ in the first and $(j-1)q+i$ in the second vector. Thus we need to find the sign of the permutation of $pq$ elements
$$
\tau: (i-1)p+j \mapsto (j-1)q+i;\quad i=1,\dots,q,\ j=1,\dots,p.
$$
Let us find the number of inversions in $\tau$: this is the number of pairs $(i_1,j_1)$ and $(i_2,j_2)$ satisfying
$$
(i_1-1)p+j_1>(i_2-1)p+j_2,\quad (j_1-1)q+i_1<(j_2-1)q+i_2
$$
or, equivalently,
\begin{equation} \label{perm}
(i_1-i_2)p>j_2-j_1,\quad (j_2-j_1)q>i_1-i_2.
\end{equation}
It follows from (\ref{perm}) that $(i_1-i_2)pq>i_1-i_2$ and hence
\begin{equation} \label{perm+}
i_1>i_2,\quad  j_2>j_1.
\end{equation}
Since $j\leq p$ and $i\leq q$, inequalities (\ref{perm+}) imply (\ref{perm}). The number of pairs satisfying (\ref{perm+}) equals 
$$
{p\choose 2} {q\choose 2} \equiv \frac{pq}{4}\ {\rm mod}\ 2.
$$
\proofend

Let us now set $p=2n, q=2$. As before, let $x$ and $y$ be the Euler classes of  the tautological bundles, and denote by $u\in H_{2n}(G_{2n}^+ (2n+2))$ and $v \in H_{4n-2}(G_{2n}^+ (2n+2))$ the respective homological Euler classes. Set $w = v \cap \dots \cap v \in H_{2n}(G_{2n}^+ (2n+2))$. It follows from Lemma \ref{multable} that 
$$
u \cap u =2, w \cap w =2 (-1)^n, u \cap w =0.
$$
As a consequence, $u$ and $w$ are linearly independent.

It is known that $H_{2n}(G_{2n}^+ (2n+2))=\Z^2$, see, e.g., \cite{Lai1}. It follows that 
one may take $u$ and $w$ for a basis in $H_{2n}(G_{2n}^+ (2n+2);\Q)$.

\begin{lemma} \label{Gss}
One has: $F_*[M]= (\chi/2) u$.
\end{lemma}

\noindent  {\bf Proof}. Note that the  induced bundles $F^*(\xi)$ and $F^*(\nu)$ are the tangent and the normal  bundles of $M$. One has:
$$
F_* ([M]) \cap u = F_* ([M]) \cdot x = [M] \cdot F^*(x) = [M] \cdot F^*(e(\xi))=
[M] \cdot e(TM) =  \chi.
$$
It is well known  that the Euler class of the normal bundle $\nu(M)$ of an embedded manifold $M$ vanishes (see, e.g., \cite{C-S}). Hence
$$
F_* ([M]) \cap w = F_* ([M]) \cdot y^n = [M] \cdot F^*(y^n) = [M] \cdot e(\nu(M))^n =  0.
$$
One can write: $F_*[M]=au+bw$ with $a,b\in \Q$, and Lemma \ref{multable} implies that $a=\chi/2$ and $b=0$.
\proofend

It follows that
\begin{equation} \label{comput}
\sigma_* F_*[M]  \cap  F_*[M] = -\left(\frac{\chi}{2}\right) u \cap \left(\frac{\chi}{2}\right) u = -\frac{\chi^2}{2}.
\end{equation}
If $\chi \neq 0$ then $F(M)$ and $\sigma F(M)$ must intersect, therefore $M$ has negatively parallel tangent spaces. This proves the first statement of Theorem \ref{nonexist}.

Let us now prove Theorem \ref{nonexist2}. 

\begin{lemma} \label{norm}
One has: $e(\nu(M)) =- 2 d$.
\end{lemma}

\noindent  {\bf Proof}.
Since $\R^4$ is contractible, one has: $[M] \cap [M]=0$. On the other hand, this homological self-intersection can be computed as follows. Choose a generic section  $\gamma$ of the normal bundle $\nu(M)$ and let $M_{\varepsilon}$ be the result of pushing $M$ slightly along this section.
Then every double point of $M$ contributes two points (with the same sign) to $M \cap M_{\varepsilon}$, and each zero of $\gamma$ contributes one point to this intersection. It  follows that $e(\nu(M)) + 2 d=0$.
\proofend

 Arguing as in Lemma \ref{Gss}, one finds that
$F_*([M])=(1-g)u+dw$, and therefore
$$
\sigma_* F_*[M]  \cap  F_*[M] = 2(d^2-(1-g)^2).
$$
This implies Theorem \ref{nonexist2}.

Let us return to the general position statement of Theorem \ref{nonexist}. If the Gauss image 
$F(M^{2n}) \subset G_{2n}^+ (2n+2)$ is an immersed submanifold  which intersects $\sigma (F(M^{2n}))$ transversally then their homological intersection equals the algebraic number of the intersection points. In particular, there are no fewer intersection points than the absolute value of the homological intersection number, see (\ref{comput}). 

Thus we need the following lemma in which ``generic"  is understood as belonging to an open and dense subset in the space of smooth maps with an appropriate topology; see, e.g., \cite{G-G}.

\begin{lemma} \label{genpos}
For a generic immersion $M^n \to  \R^{n+2}$, the tangent Gauss map $G: M \to G_n (n+2)$ is an immersion with transverse self-intersections. 
\end{lemma}

\noindent  {\bf Proof}.
Locally, $M$ is represented as the graph of a smooth map $U^n \to \R^2$ where $U$ is a domain in $\R^n$. Thus one has two functions, say, $u$ and $v$, of variables $x=(x_1,\dots,x_n)$. The tangent space to $M$ at the point corresponding to $x\in U$ is parallel to the graph of the linear map $\R^n \to \R^2$ with the $n\times 2$ matrix $A(x)=(u_{x_i}, v_{x_i}),\ i=1,\dots,n$. Hence, in an appropriate local chart of the Grassmanian, the Gauss map $G: U\to G_n (n+2)$ is given by $x\mapsto A(x)$.

The Gauss map is an immersion if  the $2n\times n$ matrix $B=(u_{x_i x_j}, v_{x_i x_j}),\break i,j=1,\dots,n$ has full rank $n$. This matrix is formed by two symmetric $n\times n$ matrices of second partial derivatives. Consider the $n(n+1)$-dimensional space of $2n\times n$ matrices formed by two symmetric $n\times n$ matrices, and let $\Delta$ be its algebraic subvariety that consists of matrices of rank $n-1$ or less. 
\smallskip

\noindent {\bf Claim}:  codim $\Delta \geq n+1$.
\smallskip

Assuming this claim, one proceeds as follows. In the space of 2-jets of smooth maps $U^n \to \R^2$, consider the algebraic subvariety $\Sigma$ of the maps for which the matrix $B$ is not of full rank. By the above claim, this subvariety has codimension $n+1$. Then the Thom transversality theorem (see, e.g., \cite{G-G}) implies that the 2-jet extension of a generic map $U^n \to \R^2$ avoids $\Sigma$, and Lemma \ref{genpos} follows.

It remain to prove the claim.
Consider a $2n\times n$ matrix $C=(S,T)$ where $S$ and $T$ are symmetric $n\times n$ matrices. By choosing an appropriate basis in $\R^n$, one can diagonalize both matrices. Assume that rank $C=n-1$. 

Assume that the first $k \leq n-1$ diagonal entries of $S$ are non-zero, and the remaining $n-k$ are equal to zero. Since the rank of $C$ is $n-1$, there are $n-k-1$ non-zero entries among  $t_{k+1,k+1},\dots, t_{n,n}$; assume that $t_{i,i} \neq 0$ for $i=k+2,\dots,n$.

Every matrix, close to $C$, can be written as $C_{\varepsilon}=(S+\varepsilon U,T+\varepsilon V)$ where $U$ and $V$ are symmetric. Then $u_{ij}$ and $v_{ij}$ are local coordinates in a neighborhood of $C$. For small $\varepsilon$, the  rank of $C_{\varepsilon}$ is not less than $n-1$. Denote the columns of $C_{\varepsilon}$ by $\eta_1,\dots,\eta_{2n}$.

If rank $C_{\varepsilon}= n-1$ then, for every $j=k+1,k+2,\dots,n+k+1$, the rank of the $n$-tuple of $n$-dimensional vectors $\eta_1,\dots,\eta_k, \eta_j, \eta_{n+k+2},\dots, \eta_{2n}$ equals $n-1$. Hence the $k+1$-st components of the vectors $\eta_j$ vanish for $j=k+1,\dots,n+k+1$, that is, $u_{k+1,j}=0$ for $j=k+1,\dots,n$ and $v_{k+1,j}=0$ for $j=1,\dots,k+1$. This gives $n+1$ equations on the variables $u_{ij}$ and $v_{ij}$, and these equations are independent of each other. This completes the proof.
\proofend

\begin{remark} \label{dim4}
{\rm  In the case $n=1$, the Grassmanian is the product of spheres:  $G_{2}^+ (4) = S^2\times S^2$. Denote the two factors by $S_1$ and $S_2$. The main result of  \cite{Bl, C-S} is that $F_*[M]=(\chi(M)/2)([S_1]+[S_2])$. The relation with classes $u$ and $v$ is as follows: $u=[S_1]+[S_2], v=[S_2]-[S_1]$. The map $d$ that takes a plane to its orthogonal complement is the identity on $S_1$ and  the antipodal involution on $S_2$. Under the identification $\R^4=\C^2$, the the space  of complex lines in $\C^2$ identifies with $\{N,S\} \times S_2$ where $N$ and $S$ are the poles of $S_1$.
}
\end{remark}

\begin{remark} \label{totreal}
{\rm  Suppose that an immersed surface $M$ in the statement of Theorem \ref{nonexist2} is totally real, that is, the tangent plane to $M$ is never a complex line. Then $F_*[M]=a [S_2]$ for some $a\in \Z$, and the algebraic number of negatively parallel tangent planes is zero. Alternatively, the normal bundle of a totally real immersed surface is isomorphic to the tangent bundle, which implies that $d^2=(1-g)^2$.

Likewise, this number is zero for an immersed surface $M$ with non-vanishing normal curvature, see \cite{Lit}. Indeed, it is proved in \cite{Lit} that the normal Euler class of $M$ is $\pm 2 \chi(M)$, that is, $d=\pm \chi(M)$. 
}
\end{remark}

\begin{remark} \label{even}
{\rm It follows from the above proof that $\chi$ is even. That this is indeed the case is a theorem of Seifert: a closed orientable embedded manifold $M^{2n}\subset \R^{2n+2}=\C^{n+1}$ has an even Euler characteristic. It is proved in \cite{Lai1} that $\chi(M)$  is twice the algebraic number of positively oriented complex points of $M$.
}
\end{remark}

\begin{remark} \label{symp}
{\rm The Grassmanian $G_{2n}^+ (2n+2) = G_{2}^+ (2n+2)$ is the space of oriented great circles in the unit sphere $S^{2n+1}$. If the space of oriented geodesics of  a  Riemannian manifold is a smooth manifold then it has a canonical symplectic structure, see, e.g.,  \cite{A-G}. Thus $G_{2n}^+ (2n+2)$ is a symplectic manifold, and $\omega^{2n}\neq 0$ where $\omega$ is the symplectic form. In the previous notation, the cohomology class of $\omega$ is $y$.
}
\end{remark}

\begin{remark} \label{inter}
{\rm  Let $M^{m}\subset \R^{m+2}$ be a closed immersed manifold and $e\in \R^{m+2}$ an arbitrary non-zero vector. Then the Gauss image $F(M)  \subset G_m^+ (m+2)$ intersects $N_e$, and for a generic $e$, in at least rk $H_*(M)$ points. Indeed, the tangent space $T_x M$ is orthogonal to $e$ if and only if the restriction of the ``height" function $f(y)=y \cdot e$ on $M$ has a critical point at $x$. In contrast, the Gauss image $F(M)$ can be disjoint from the set $P_f\subset G_m^+ (m+2)$: for example, this is the case when $M$ lies in a hyperplane and the vector $f$ is orthogonal to this hyperplane.
}
\end{remark}

\section{Examples of skew branes} \label{pf2}

\subsection{Odd-dimensional skew sphere}\label{sph}

We  construct an odd-dimensional skew sphere  described in Theorem \ref{exist}.

Let $M^{m-1} \subset \R^{m}$ be a smooth strictly convex closed
hypersurface containing the origin. Denote by $h: S^{m-1} \to \R$ the support 
function of $M$, that is, $h(x)$ is the distance from the origin to the tangent hyperplane to $M$ for which the unit vector $x \in S^{m-1}$ is the outward normal.  Extend $h$ to $\R^{m}$ as a homogeneous function of degree 1. By Euler formula, the latter implies that 
\begin{equation} \label{Euler}
x h_x=h,\quad   x h_{xx}=0.
\end{equation}
Here and elsewhere we use following convention: for vectors $x,y$, a function $g$,
one has:
$$
x y=\sum_{i=1}^m x_i y_i,\ \ x dy = \sum_{i=1}^m x_i dy_i,\ \ g_x=(g_{x_1},\dots,g_{x_m}),
$$
$$
x g_{xx}=(\sum_{i=1}^m x_i g_{x_i x_1},\dots, \sum_{i=1}^m x_i g_{x_i x_m}),
$$
etc. 

Thus the hypersurface $M$ is parameterized by the unit sphere as follows.

\begin{lemma} \label{param}
Let $y(x) \in M$ be the point at  which the outward unit normal to $M$ is $x\in S^{m-1}$. Then  $y= h_x$.
\end{lemma} 

\noindent  {\bf Proof}. 
Consider the hypersurface  $M'$ given by the formula $y(x)= h_x$. To prove that the tangent hyperplane to $M'$ at $y$ is orthogonal to $x$, one needs to check  that the 1-form $x dy$ vanishes on $M'$. Indeed, $xy=h$ by (\ref{Euler}), and hence
$$
x dy= d(xy)-y dx = dh-h_x dx=0.
$$
Thus $h$ is the support function of $M'$ and therefore $M'=M$. 
\proofend

It follows from (\ref{Euler}) that, at point $x \in S^{m-1}$,   the linear operator $h_{xx}$ annihilates  the normal direction $x$ and preserves the tangent space $T_x S^{m-1}$. The convexity of $M$ implies that, on the tangent space, $h_{xx}$ is non-degenerate (for example, $h_{xx}$ is the identity if $h=1$ and hence $M$ is the unit sphere). Denote the restriction of $h^{-1}_{xx}$ on $T_x S^{m-1}$ by $A(x)$ and extend $A$ to the normal direction, spanned  by $x$, as the zero map.

We will construct  a skew brane in $\R^{m+1}$ as a section of the vertical cylinder over $M\subset \R^m \subset \R^{m+1}$. More specifically, let $f: S^{m-1} \to {\R}$ be a smooth function. Define
$$
N=\{(h_x (x),f(x)) |\ x\in S^{m-1} \} \subset \R^m \times \R = \R^{m+1},
$$
and let $\phi: S^{m-1} \to N$ be the parameterization map.

Let us describe when the tangent spaces to $N$ are parallel. Assume that $f$ is extended to $\R^m$ as a homogeneous function of some degree $k$ (whose value is of no importance).

\begin{lemma} \label{parall}
The tangent spaces to $N$ are parallel at distinct points  $\phi(x_1)$ and $\phi(x_2)$ if and only if $x_2=-x_1$ and $A(f_x) (x_2)=A(f_x) (x_1)$.
\end{lemma} 

\noindent  {\bf Proof}. Let us describe the normal 2-plane to $N$ at  point $\phi(x)$. This plane is generated by the vector $(x,0)$ and a vector $(\xi,1)$ where $\xi \in T_x S^{m-1}$. We claim that $\xi =-A(f_x)$.

Indeed, let $v\in T_x S^{m-1}$ be a test vector, $v=\sum_j v_j \partial_{x_j}$. Then 
$$
d\phi(v)=\sum_j (v_j h_{x_i x_j},v_j f_{x_j}),\ \ i=1,...,m.
$$
It follows that $d\phi(v)$ is orthogonal to $(\xi,1)$ if and only if $v h_{xx} \xi + vf_x=0,$
or $v (h_{xx} \xi + f_x)=0.$ Since $v$ is an arbitrary tangent vector to $S^{m-1}$, the projection of the vector $h_{xx} \xi + f_x$ to $T_x S^{m-1}$ is zero. This implies the claim.

Finally, the span of  vectors $(\xi(x_2),1)$ and $(x_2,0)$ coincides with that of $(\xi(x_1),1)$ and $(x_1,0)$ if and only if $x_2=-x_1$ and  $\xi(x_2)=\xi(x_1)+tx_1$.
Since $\xi(x_2)$ and $\xi(x_1)$ are orthogonal to $x_1$, one has $\xi(x_2)=\xi(x_1)$.
\proofend

Thus  $N$ is a skew brane if and only if
\begin{equation} \label{grad1}
A(-x) (f_x (-x)) \neq A(x) ( f_x(x)) 
\end{equation}
for all $x\in S^{m-1}$.
Decompose the functions $h$ and $f$ into the even and odd components with respect to the antipodal involution of the sphere $x\mapsto -x$:
$$
h=h_{ev}+h_{odd},\quad f=f_{ev}+f_{odd}.
$$
Then one has a similar decompositions of the operator $A$ and the gradient
vector field $f_x$. Note that
$$
(f_x)_{ev}=(f_{odd})_x,\ (f_x)_{odd}=(f_{ev})_x.
$$
Decomposing  (\ref{grad1}) into even and odd parts yields:
\begin{equation} \label{grad2}
A_{ev} ((f_{ev})_x) + A_{odd}  ((f_{odd})_x) \neq 0 
\end{equation}
for all $x\in S^{m-1}$. 

If $M$ is sufficiently close  to the sphere 
then $h$ is close to 1 and $A$ is close to the identity. Thus we may assume that $A_{ev}$  is invertible and rewrite (\ref{grad2}) as
\begin{equation} \label{grad3}
(f_{ev})_x + B((f_{odd})_x) \neq 0  
\end{equation}
 where $B(x)=A_{ev}^{-1}A_{odd}$ is an odd field of linear maps of the tangent spaces $T_x S^{m-1}$. 

Note that (\ref{grad3}) cannot hold for all $x$ if $h$, and therefore $A$, is 
even: then $B=0$ and the function $f_{ev}$ must have critical points on the sphere. Note also that (\ref{grad3}) cannot hold if $m-1$ is even since there exist no non-vanishing vector fields on even-dimensional spheres. From now on, assume that $m=2n$.

\begin{lemma} \label{max}
Let $U$ be a contractible domain in $\R^m$ and $B$ a smooth field of linear maps of the cotangent spaces to $U$ satisfying the following property: for every smooth function $f:U\to \R$ there exists a smooth function $g:U\to \R$ such that $B(df)=dg$. Then $B=c\ {\rm Id}$ for some constant $c$.
\end{lemma}

\noindent  {\bf Proof}. In local coordinates $x_1,\dots,x_m$, the linear map $B$ is given by a matrix $b_{ij}(x)$. The condition is that, for every $f$, the 1-form $\sum_{ij} b_{ij} f_{x_j} dx_i$ is closed.

Choose $j\in\{1,\dots,m\}$ and let $f(x)=x_j$. Then the 1-form $\sum_i b_{ij}(x) dx_i$ is closed, and hence there exists a function $g^j (x)$ such that $b_{ij}=g^j_{x_i}$. The 1-form $\sum_{ij} g^j_{x_i} f_{x_j} dx_i = \sum_j f_{x_j} d g^j$ is closed for every $f$, therefore so is $\sum_j g^j d (f_{x_j})$.

Again fix $j\in\{1,\dots,m\}$ and let $f(x)=x_j^2/2$. Then the 1-form $g^j d x_j$ is closed, hence $g^j$ depends only on $x_j$. Finally, fix  $j,k\in\{1,\dots,m\}$, $j\neq k$, and let $f(x)=x_j x_k$. Then the 1-form $g^j (x_j) d x_k + g^k (x_k) dx_j$ is closed, and therefore $g^j_{x_j} = g^k_{x_k}=c$. It follows that $b_{ij}=c \delta_{ij}$, as claimed.
\proofend

One has the following corollary.

\begin{corollary} \label{nonexact}
If $M$ is not centrally symmetric then there exists a function $f_{odd}$ on $S^{m-1}$ such that  $B((f_{odd})_x)$ is not a gradient vector field.
\end{corollary}

\noindent  {\bf Proof}. A parallel translation of the origin changes the support function $h$ by addition of  a linear function and does not affect $h_{xx}$. The hypersurface $M$ is centrally symmetric if and only if, after a parallel translation,   $h$ is even, that is, $A$ is even. Since $M$ is not centrally symmetric, $A_{odd} \neq 0$ and $B \neq 0$. 

Extend $f$ to $\R^m$ as a homogeneous  function of degree $0$ and  identify the tangent and cotangent spaces by the Euclidean structure. Let $U$ be a small domain on $S^{m-1}$. Then every function $f: U\to \R$ can be extended to $S^{m-1}$ as an odd function.  If $B(f_x)$ is a gradient for every such $f$ then, by Lemma \ref{max}, $B=c\ {\rm Id}$ in $U$. This equality then must hold everywhere on the sphere. But $B$  is odd, hence $B=0$, a contradiction. 
\proofend

We make the following conjecture.

\begin{conjecture} \label{cylinder}
For every non-centrally symmetric $M^{2n-1}\subset \R^{2n}$ there exists a skew brane $N^{2n-1}\subset \R^{2n+1}=\R^{2n}\times \R$ which is a section of the vertical cylinder over $M$.
\end{conjecture}

Conjecture \ref{cylinder}  follows from a more general conjecture of independent interest. 

\begin{conjecture} \label{oneform}
Let $M$ be a closed manifold with zero Euler  characteristic, $\alpha$ a non-closed differential 1-form on it. Then there exists a smooth function $f:M\to {\bf R}$ such that the 1-form $df+\alpha$ has no zeroes on $M$.
\end{conjecture}
 
\noindent  {\bf Proof of Conjecture \ref{cylinder} from Conjecture \ref{oneform}}.
One may rephrase (\ref{grad3}) as 
\begin{equation} \label{reph}
d (f_{ev}) + \alpha \neq 0
\end{equation}
everywhere on $S^{2n-1}$; here $\alpha$ is  the 1-form dual to the vector field $B((f_{odd})_x)$. By Corollary \ref{nonexact}, we may choose $f_{odd}: S^{2n-1}\to \R$ so that $\alpha$ is not closed. Consider the quotient space $M=\RP^{2n-1}$. Then $\alpha$ descends as a non-closed 1-form $\bar \alpha$ on $M$. Conjecture \ref{oneform} implies that there exists a function $\bar f: \RP^{2n-1} \to \R$ such that 
$d \bar f+\bar \alpha$ is nowhere vanishing. Then $\bar f$ lifts to an even function on $S^{2n-1}$ for which (\ref{reph}) holds.
\proofend

\begin{remark} \label{cyldim1}
{\rm In dimension one, Conjecture \ref{oneform} holds. Indeed, $\alpha = g(x) dx$ with $c:=\int_0^1 g(x) dx \neq 0$; here $x$ is the coordinate on the circle $\R/\Z$. Then $g-c$ is a derivative, $g(x)-c=-f'(x)$, and $df+\alpha=cdx$. 
It follows that, in dimension one, Conjecture \ref{cylinder} holds as well: this the Cylinder Lemma of \cite{G-S}.
}
\end{remark}

At present, we cannot prove Conjecture \ref{oneform}. To prove Theorem \ref{exist}, we will construct functions $h$ and $f$ on $S^{2n-1}$ so that (\ref{grad3}) holds everywhere on $S^{2n-1}$. We set $h=1+\varepsilon g$ where $g$ is an odd function and $\varepsilon$ is a small positive real. Then  $h_{xx}=Id+ \varepsilon g_{xx} + O(\varepsilon^2)$ and $A=Id - \varepsilon g_{xx}+ O(\varepsilon^2)$. It follows that $B=
-\varepsilon g_{xx} + O(\varepsilon^2)$, and (\ref{grad3}) can be rewritten as
\begin{equation} \label{grad4}
(f_{ev})_x \neq \varepsilon g_{xx} ((f_{odd})_x) + O(\varepsilon^2).
\end{equation}
Our strategy is to construct a function $f$ and an odd function $g$ so that, for sufficiently small  $\varepsilon$ and a certain constant $c>0$, 
\begin{equation} \label{grad5}
|(f_{ev})_x - \varepsilon g_{xx} ((f_{odd})_x)| > c \varepsilon
\end{equation}
everywhere on the sphere. This would imply (\ref{grad4}).

Let $(z_1,\dots,z_n)$ be coordinates in $\C^n$, and let  $f_{ev}=\sum_i a_i |z_i|^2$
with generic  real  $a_i$.  Denote by $\xi$ the unit Hopf vector field on $S^{2n-1}$.
Then $f_{ev}$ is an even, $\xi$-invariant Morse-Bott function on $S^{2n-1}$ with $n$ critical Hopf circles $C_1,\dots,C_n$. 

\begin{lemma} \label{other}
There exist odd smooth functions $f_{odd}$ and $g$ on  $S^{2n-1}$ such that $g_{xx} ((f_{odd})_x) \cdot \xi =2$ on each critical circle $C_i, i=1,\dots,n$.
\end{lemma}

Assuming this lemma, denote the vector field $g_{xx} ((f_{odd})_x)$ by $v$. Let $c$ be a constant such that $|v(x)|<c$ everywhere on the sphere. Given $\varepsilon>0$, let $U_{\varepsilon}$ be a neighborhood of the critical set $C_1 \cup \dots \cup C_n$ such that  $|(f_{ev})_x| > 2 \varepsilon c$ outside of $U_{\varepsilon}$. We claim that (\ref{grad5}) holds outside of $U_{\varepsilon}$. Indeed, $|(f_{ev})_x -\varepsilon v| > 2 \varepsilon c - \varepsilon c=\varepsilon c$, as claimed.

Next, consider the situation inside $U_{\varepsilon}$. Since $v \cdot \xi =2$ on the critical circles, one has $v \cdot \xi > 1$ inside $U_{\varepsilon}$ for sufficiently small $\varepsilon$. Note that $(f_{ev})_x$ is orthogonal to $\xi$. Therefore
$|((f_{ev})_x - \varepsilon v) \cdot \xi| = \varepsilon v \cdot \xi > \varepsilon$, and hence $|((f_{ev})_x - \varepsilon v)| > \varepsilon$ inside $U_{\varepsilon}$. In particular, (\ref{grad5}) holds. 

Thus,  Theorem \ref{exist} will follow, once we prove Lemma \ref{other}.
\medskip

\noindent  {\bf Proof of Lemma \ref{other}}. We construct the desired functions in a neighborhood of each critical circle and then extend them to the sphere. Let $C$ be one such Hopf circle. Consider $\C^2\subset \C^n$ such that $C$ is a Hopf circle therein. Choose coordinates $(z_1,z_2)$  in $\C^2$ so that $z_i = x_i +\sqrt{-1}y_i,\ i=1,2$, and $C$ is given by $x_2=y_2=0$. The Hopf field is given by the formula:
$$
\xi=-y_1 \partial x_1 + x_1 \partial y_1 - y_2 \partial x_2 + x_2 \partial y_2.
$$
The functions $f=f_{odd}$ and $g$ will depend only on $(x_1,y_1,x_2,y_2)$.

It is straightforward to compute the operator $g_{xx}$, the gradient $f_x$ and the dot product $g_{xx} (f_x) \cdot \xi$. The answer is as follows. Introduce the differential operator depending on $f$:
$$
D= f_{x_1} \partial x_1 + f_{y_1} \partial y_1 + f_{x_2} \partial x_2+ f_{y_2} \partial y_2.
$$
Then, 
\begin{equation} \label{oper}
g_{xx} (f_x) \cdot \xi = x_1 D(g_{y_1}) -y_1 D(g_{x_1}) + x_2 D(g_{y_2}) -y_2 D(g_{x_2}).
\end{equation}
Let us look for $f$ and $g$ in the following form:
$$
f=a(x_1,y_1)+x_2 b(x_1,y_1)+y_2 d(x_1,y_1),\  g=u(x_1,y_1) + x_2 v(x_1,y_1) + y_2 w(x_1,y_1)
$$
where $a,u$ are odd and homogeneous of degree 1, $b,d,v,w$ are even and homogeneous of degree 0. Switch to polar coordinates $x_1=r\cos \alpha, y_1=r\sin \alpha$. On the circle $C$, one has $x_2=y_2=0, r=1$, and $\alpha$ is a coordinate. Then (\ref{oper}) becomes
\begin{equation} \label{oper1}
bv'+dw'+(u+u'')a'
\end{equation}
where prime is $d/d \alpha$. Set: 
$$
a=u=0,\ \ b=-w=\cos 2 \alpha,\ \ d=v=\sin 2 \alpha,
$$
and the expression (\ref{oper1}) gets identically equal to $2$. In terms of the Cartesian coordinates,
$$
f= x_2(x_1^2-y_1^2)-2y_2x_1y_1,\ \  
g=2x_2x_1y_1-y_2(x_1^2-y_1^2)
$$
on the unit sphere. This completes the construction.
\proofend

\subsection{Skew torus in four-dimensional space} \label{tor}

Start with the standard torus $T^2_0$ given by a parameterization 
$$
(\alpha,\beta) \mapsto (\cos \alpha, \sin \alpha, \cos \beta, \sin \beta);
$$
this torus is an orbit of a 1-parameter group of isometries of $\R^4$.  One has an action of the group $\Z_2^2$ on $T^2_0$ generated by 
\begin{equation} \label{group}
S_1: (\alpha,\beta) \mapsto (\alpha+\pi,\beta),\ \ S_2:(\alpha,\beta) \mapsto (\alpha,\beta+\pi).
\end{equation}
 The tangent planes to $T^2_0$ at points $(\alpha_1,\beta_1)$ and $(\alpha_2,\beta_2)$ are parallel if and only if  $(\alpha_1,\beta_1)$ and $(\alpha_2,\beta_2)$ belong to the same $\Z_2^2$-orbit.  The Gauss map $T^2_0 \to G_2(4)$ is an immersion, namely, a four-fold covering of its image. 

Consider  a new torus $T^2$, perturbation of $T^2_0$, given by the formula
$$
(\alpha,\beta) \mapsto (u \cos \alpha, u \sin \alpha, v \cos \beta, v \sin \beta)
$$
where $u=1+\varepsilon f(\alpha,\beta), v=1+\varepsilon g(\alpha,\beta)$. We will construct functions $f$ and $g$ such that, for $\varepsilon$ small enough, $T^2$ has no pairs of parallel tangent planes.

For a sufficiently small $\varepsilon$, the Gauss map for $T^2$ is an immersion.  Given $(\alpha,\beta) \in T^2$, assume that the tangent plane at point $(\bar \alpha,\bar \beta)$ is parallel to that at $(\alpha,\beta)$. Then $(\bar \alpha,\bar \beta)$ is close to either of the three points: $S_1(\alpha,\beta), S_2(\alpha,\beta)$ or $S_1 S_2(\alpha,\beta)$. 
Consider the linearization of the equations that express the fact that $T_{(\alpha,\beta)} T^2$ and $T_{(\bar \alpha,\bar \beta)} T^2$ are parallel: this linearization is obtained by ignoring the terms of order 2 and higher in $\varepsilon$. 

\begin{lemma} \label{paral}
In the linear approximation in $\varepsilon$, if the planes $T_{(\alpha,\beta)} T^2$ and $T_{(\bar \alpha,\bar \beta)} T^2$ are parallel then one of the following three systems of two equations holds for some $(\alpha,\beta)$:
\begin{equation} \label{case1}
f_{\beta}(\alpha+\pi, \beta)+f_{\beta}(\alpha,\beta) =0,\ \ g_{\alpha}(\alpha+\pi, \beta)+g_{\alpha}(\alpha,\beta) =0;
\end{equation}
\begin{equation} \label{case2}
f_{\beta}(\alpha, \beta+\pi)+f_{\beta}(\alpha,\beta) =0,\ \ g_{\alpha}(\alpha, \beta)+g_{\alpha}(\alpha,\beta+\pi) =0;
\end{equation}
\begin{equation} \label{case3}
f_{\beta}(\alpha+\pi, \beta+\pi)-f_{\beta}(\alpha,\beta) =0,\ \ g_{\alpha}(\alpha+\pi, \beta+\pi)-g_{\alpha}(\alpha,\beta) =0.
\end{equation}
\end{lemma} 

\noindent  {\bf Proof}. Choose a basis $e_1,e_2,e_3,e_4$ in $\R^4$. Given a plane  $E\subset \R^4$, choose a  frame $(f_1,f_2)$ in $E$ and consider the bivector $f_1\wedge f_2\in \Lambda^2(\R^4)$. This bivector is uniquely defined, up to a factor, by $E$. One can write
$$
f_1\wedge f_2 = \sum_{1\leq i<j\leq 4} p_{ij} e_i\wedge e_j;
$$
$p_{ij}$ are called the Plucker coordinates of the plane $E$. The Plucker coordinates are
defined up to a common factor and are not independent: they satisfy the identity 
$p_{12}p_{34}-p_{13}p_{24}+p_{14}p_{23}=0$. Reversing orientation of the plane changes the signs of all Plucker coordinates.

It is convenient to change coordinates:
$$
x_1=p_{12}+p_{34},\ x_2=p_{23}+p_{14},\ x_3=-p_{13}+p_{24},
$$
$$
y_1=p_{12}-p_{34},\ y_2=p_{23}-p_{14},\ y_3=-p_{13}-p_{24}.
$$
Then the tangent plane $T_{(\alpha,\beta)} T^2$ has the following Plucker coordinates (we continue ignoring the terms of order $\varepsilon^2$ and higher):
$$x_1=\varepsilon(g_\alpha(\alpha,\beta)-f_\beta(\alpha,\beta)),\ x_2=-\sin (\alpha+\beta) + O(\varepsilon),\ x_3=\cos (\alpha+\beta) + O(\varepsilon),
$$
$$y_1=-\varepsilon(g_\alpha(\alpha,\beta)+f_\beta(\alpha,\beta)),\ y_2=\sin (\alpha-\beta) + O(\varepsilon),\ y_3=-\cos (\alpha-\beta) + O(\varepsilon).
$$

Assume that $T_{(\alpha,\beta)} T^2$ and $T_{(\bar \alpha,\bar \beta)} T^2$ are parallel.   The three systems of the lemma correspond to the following three cases:
$$
(\bar \alpha,\bar \beta)=S_1(\alpha,\beta)+ O(\varepsilon),\ (\bar \alpha,\bar \beta)=S_2(\alpha,\beta)+ O(\varepsilon),\ (\bar \alpha,\bar \beta)=S_1 S_2(\alpha,\beta)+ O(\varepsilon).
$$ 
The three cases  being similar, consider the first: 
$$
\bar \alpha=\alpha +\pi + O(\varepsilon),\ \bar \beta=\beta + O(\varepsilon).
$$ 
Denote the Plucker coordinates of the plane $T_{(\bar \alpha,\bar \beta)} T^2$ by $\bar x_i$ and $\bar y_i$, $i=1,2,3$. Then the vectors $(\bar x_i, \bar y_i)$ and $(x_i, y_i)$ are proportional. Since the zero-order terms in $\varepsilon$ are opposite, one has
$$
\bar x_i = (-1+\varepsilon c) x_i + O(\varepsilon^2),\ \bar y_i = (-1+\varepsilon c) y_i + O(\varepsilon^2)
$$
for some real $c$. In particular, for $i=1$, this implies:
$$
g_\alpha(\alpha+\pi,\ \beta)-f_\beta (\alpha+\pi,\beta)=f_\beta(\alpha, \beta)-g_\alpha(\alpha, \beta),
$$
$$
-g_\alpha(\alpha+\pi,\beta)-f_\beta(\alpha+\pi,\beta)=f_\beta(\alpha, \beta)+g_\alpha(\alpha, \beta),
$$
and (\ref{case1}) follows.
\proofend

The number of solutions of the linearized system provides an upper bound on the number of genuine solutions; in particular, if the linearized system does not have solutions then neither does the original system. This principle implies that
it suffices to find functions $f$ and $g$ on the torus for which systems of equations (\ref{case1}), (\ref{case2}) and (\ref{case3}) do not hold for all $(\alpha, \beta)$. 

More specifically, here and in Section \ref{immsph}, we use the following lemma.

\begin{lemma} \label{just}
Let $f,g$ and $h$ be continuous functions on a compact space $M$ such that 
for all $x\in M$ and all sufficiently small  $\varepsilon \geq 0$ one has: $f(x)+\varepsilon g(x) >0$. Then $f(x)+\varepsilon g(x) + \varepsilon^2 h(x) >0$
for all $x\in M$ and all sufficiently small  $\varepsilon \geq 0$.
\end{lemma}

\noindent  {\bf Proof}. It suffices to prove that there exists a constant $c>0$ such that 
\begin{equation} \label{linest}
f(x)+\varepsilon g(x) >c \varepsilon
\end{equation}
 for all $x\in M$ and all sufficiently small  $\varepsilon \geq 0$. Indeed, if (\ref{linest}) holds then 
$$
f(x)+\varepsilon g(x) + \varepsilon^2 h(x) > c \varepsilon - \varepsilon^2 |h(x)| >0
$$
for $\varepsilon < c/b$ where $b={\rm max}_{x\in M} |h(x)|$.

Now we prove (\ref{linest}). Let $N \subset M$ be the zero locus of $f(x)$. Then $N$ is  compact. By assumption, the restriction of $g$ on $N$ is everywhere positive. Set:  
 $0<2a={\rm min}_{x\in N} g(x)$ and $C={\rm max}_{x\in M} |g(x)|$. For a given $\varepsilon$, let $U_{\varepsilon}$ be a neighborhood of $N$ such that $f(x) > 2 \varepsilon C$ for all $x$ outside of $U_{\varepsilon}$. Then, outside of $U_{\varepsilon}$, one has: 
$ f(x)+\varepsilon g(x) \geq f(x) - \varepsilon C >  \varepsilon C. $

Consider the situation inside $U_{\varepsilon}$. If $\varepsilon$ is small enough then $g(x) > a$ for all $x \in U_{\varepsilon}$. It follows that
$f(x)+\varepsilon g(x) \geq \varepsilon g(x) >  \varepsilon a$,
and we are done.
\proofend

Now we are ready to finish the proof of of Theorem \ref{torus}.
Regarding the $\Z_2^2$ action (\ref{group}), every function $f(\alpha,\beta)$ can be decomposed $f=f^{0,0}+f^{1,0}+f^{0,1}+f^{1,1}$ where $f^{0,0}$ is even with respect to $S_1$ and $S_2$, $f^{1,0}$ is odd with respect to $S_1$ and even with respect to $S_2$,  $f^{0,1}$ is even with respect to $S_1$ and odd with respect to $S_2$, and $f^{1,1}$ is odd with respect to $S_1$ and $S_2$. This decomposition is preserved by partial differentiation with respect to $\alpha$ and $\beta$.

Then equations (\ref{case1}), (\ref{case2}) and (\ref{case3}) are equivalent, respectively, to
\begin{equation} \label{case11}
f_{\beta}^{0,0} + f_{\beta}^{0,1}=0,\ \ g_{\alpha}^{0,0} + g_{\alpha}^{0,1}=0;
\end{equation}
\begin{equation} \label{case12}
f_{\beta}^{0,0} + f_{\beta}^{1,0}=0,\ \ g_{\alpha}^{0,0} + g_{\alpha}^{1,0}=0;
\end{equation}
and
\begin{equation} \label{case13}
f_{\beta}^{1,0} + f_{\beta}^{0,1}=0,\ \ g_{\alpha}^{1,0} + g_{\alpha}^{0,1}=0.
\end{equation}

Set $f^{0,0} = \cos(2\alpha+2\beta), g^{0,0} = \sin(2\alpha+2\beta)$. Then 
\begin{equation} \label{iden}
(f_{\beta}^{0,0})^2+(g_{\alpha}^{0,0})^2=4
\end{equation}
for all $(\alpha, \beta)$. Set $f^{0,1}=g^{0,1}=0$ and $f^{1,0}= \cos(\alpha+2\beta), g^{1,0} = 2\sin(\alpha+2\beta)$. Then $(f_{\beta}^{1,0})^2+(g_{\alpha}^{1,0})^2=4$ identically on the torus, and hence system (\ref{case13}) does not hold for all $(\alpha, \beta)$. Finally,  multiply $f^{1,0}$ and $g^{1,0}$ by a sufficiently small constant. By continuity, it follows from  (\ref{iden}) that systems  (\ref{case11}) and  (\ref{case12}) do not hold for all $(\alpha, \beta)$. This completes the proof of Theorem \ref{torus}.

\subsection{Immersed sphere in four-dimensional space} \label{immsph}

Now we prove Theorem \ref{2sph}; similarly to the preceding section, we use Lemma \ref{just}. Consider the following perturbation of $M_0$, the surface $M$ given by the equation
$$
(\alpha,h) \mapsto (1-h^2) (\cos \alpha, \sin \alpha, h \cos \alpha - \varepsilon g(\alpha) \sin \alpha, h \sin \alpha + \varepsilon g(\alpha) \cos \alpha)
$$
where $g$ is a function to be chosen later. 

As in the case of a skew torus, we compute the Plucker coordinates of the tangent planes to $M$, in the linear approximation in $\varepsilon$. We use the notation from the proof of Lemma \ref{paral}.

\begin{lemma} \label{Plcoord}
Up to the terms of order $\varepsilon^2$ and higher, the Plucker coordinates of the oriented plane $T_{(\alpha,h)} M$ are
$$
x_1=h(1+3h^2) - \varepsilon (1-3h^2) g'(\alpha),\ x_2 = \cos (2 \alpha) ((1-h^2) + 2\varepsilon h g'(\alpha)), 
$$
$$
x_3 =  \sin (2 \alpha) ((1-h^2) + 2\varepsilon hg'(\alpha)),\ 
y_1 = 3h(1-h^2) + \varepsilon (1-3h^2) g'(\alpha), 
$$
$$
y_2 = (1-5h^2) - 2\varepsilon hg'(\alpha),\ y_3 = 4\varepsilon hg(\alpha).
$$
\end{lemma} 

The proof is a straightforward computation which we omit. In particular, we find when the tangent planes to the unperturbed sphere $M_0$ are negatively parallel.

\begin{corollary} \label{unpert}
The planes $T_{(\alpha_1,h_1)} M_0$ and $T_{(\alpha_2,h_2)} M_0$ are negatively parallel if and only if $h_1=- h_2=1/\sqrt{5}, 2\alpha_2=2\alpha_1+\pi$.
\end{corollary}

\noindent  {\bf Proof}. The planes are parallel if their Plucker coordinates are proportional. Let $\lambda$ be the proportionality factor. One has 6 equations, naturally  labeled  by the respective Plucker coordinate. Equation $(x_1)$ and  $(y_1)$ read, respectively:
\begin{equation} \label{tech1}
h_2(1+3h_2^2)=\lambda h_1(1+3h_1^2),\ h_2(1-h_2^2)=\lambda h_1(1-h_1^2),
\end{equation}
and equations $(x_2)$ and $(x_3)$ imply:
\begin{equation} \label{tech2}
(1-h_2^2)=|\lambda| (1-h_1^2).
\end{equation}
 If the planes are negative parallel then $\lambda <0$. It follows from (\ref{tech1}) and (\ref{tech2})  that $\lambda=-1$ and $h_2=-h_1$. Equation $(y_2)$ implies now that $1-5h_1^2=0$. Finally, equations $(x_2)$ and $(x_3)$ imply that $\cos (2 \alpha_2)=-\cos (2 \alpha_1)$ and $ \sin (2 \alpha_2)=-\sin (2 \alpha_1)$.
\proofend

Now we repeat this computation, taking the terms linear in $\varepsilon$ into account. Write:
$$
h_1=\frac{1}{\sqrt{5}} + \varepsilon a,\ h_2=-\frac{1}{\sqrt{5}} + \varepsilon b,\ \lambda =(-1+\varepsilon c),\ 2\alpha_2=2\alpha_1+\pi+\varepsilon d.
$$
Then Lemma \ref{Plcoord} yields 6 equations in unknowns $a,b,c,d$ and $\alpha=\alpha_1$, again labeled by the Plucker coordinates. In particular, Eq. $(y_3)$ gives:
\begin{equation} \label{col1}
g(\beta)=g(\alpha),
\end{equation}
and the linear combination 
$$
3\ {\rm Eq.} (x_1) - 7\ {\rm Eq.} (y_1) - 15 \left(\cos 2 \alpha\ {\rm Eq.} (x_2) + \sin 2 \alpha\ {\rm Eq.} (x_3) \right)
$$
yields:
\begin{equation} \label{col2}
g'(\beta)=-g'(\alpha);
\end{equation}
here $\beta = \alpha_2$ mod $\varepsilon$, that is, $2\beta=2\alpha + \pi$.

Finally, set $g(\alpha)=\sin 2 \alpha + \sin 4 \alpha$. Then equations (\ref{col1}) and (\ref{col2}) imply, respectively, $\sin 2 \alpha=0$ and $ \cos 4 \alpha =0$, that is,
are not compatible for all $\alpha$, and Theorem \ref{2sph} follows.

\end{document}